\documentclass[12pt,a4paper]{article}
\usepackage[ansinew]{inputenc}
\usepackage{amsmath}
\usepackage{amsfonts}
\usepackage{amssymb}
\usepackage{amsthm}
\usepackage [dvips] {graphicx}
\newtheorem{prop} {Proposition} 
 
\newtheorem{thm} [prop] {Theorem} 
\newtheorem{cor} [prop] {Corollary} 

\newtheorem{df}{Definition} 

\theoremstyle{remark}
\newtheorem{rmq}{Remark} 
\newtheorem{example}{Example}

\author{S. Dugowson\\\begin{small}Institut Supérieur de Mécanique de Paris\end{small}\\ \begin{small}stephane.dugowson@supmeca.fr\end{small}}

\title{On Connectivity Spaces \\ \begin{small}(to appear : 2010)\end{small}}

\date{29.08.2008}

\begin{document}

\markboth{S. Dugowson}
{On Connectivity Spaces }

\maketitle

\begin{abstract}

This paper presents some basic facts about the so-called connectivity spaces. In particular, it studies the generation of connectivity structures, the existence of limits and colimits in the main categories of connectivity spaces, the  closed monoi\-dal category structure given by the so-called tensor product on integral connectivity spaces; it defines homotopy for connectivity spaces and mention briefly related difficulties; it defines smash product of pointed integral connectivity spaces and shows that this operation results in a closed monoidal category with such spaces as objects. Then, it studies finite connectivity spaces, associating a directed acyclic graph with each such space and then defining a new numerical invariant for links: the connectivity order. Finally, it mentions the not very wellknown Brunn-Debrunner-Kanenobu theorem which asserts that every finite integral connectivity space can be represented by a link.\\



\end{abstract}

\section*{Introduction}

Connectivity spaces are kinds of topological objects which have not yet received very great attention. This paper presents results we have recently obtained relating to those spaces. In the first section, we recall their definition. The second section is about the generation of connectivity structures from a given family of subsets we wish to consider as connected. The third section is about categorical constructions in the main categories of connectivity spaces, seen as peculiar cases of the so-called categories with lattices of structures. The  fourth section studies the closed monoi\-dal category structure given by the so-called tensor product on integral connectivity spaces. The fifth section defines homotopy for connectivity spaces and briefly mentions  some difficulties related to this notion. The sixth section is devoted to pointed integral connectivity spaces and to smash product between such spaces. In the last section we study finite connectivity spaces, associating a directed acyclic graph with each such space and then defining a new numerical invariant for links: the connectivity index. Finally, it mentions the not very well-known Brunn-Debrunner-Kanenobu theorem which asserts that every finite integral connectivity space can be represented by a link in the space $\mathbf{R}^3$ (or in $\mathbf{S}^3$).

\section*{Notations}

If $X$ is a set, the set of subsets of $X$ is denoted by $\mathcal{P}(X)$ or $\mathcal{P}_X$, and the set $\mathcal{P}(\mathcal{P}_X)$ by $\mathcal{Q}_X$. 
For any $\mathcal{A}\in\mathcal{Q}_X$, $\mathcal{A}^\bullet$ denotes the set $\{A\in \mathcal{A}, \mathrm{card}(A)\geq 2\}$.
If $\sim$ is an equivalence relation on $X$, the equivalence class of $x\in X$ is denoted by $\tilde{x}$. If $Y$ is a subset of $X$, $\sim_Y$ denotes the equivalence relation defined on $X$ by $a\sim_Y b$ if and only if $a=b$ or $(a,b)\in Y^2$, and $X/Y$  denotes the quotient $X/\sim_Y$.

\section{Definitions, Examples}

Let us recall the definition of connectivity spaces and connectivity morphisms \cite{Borger:1983,Dugowson:2007c}.

\begin{df} [Connectivity spaces]  A \emph{connectivity space} is a pair \linebreak $(X,\mathcal{K})$ where $X$ is a set and $\mathcal{K}$ is a set of subsets of  $X$ such that $\emptyset\in\mathcal{K}$ and
\begin{displaymath}
\forall \mathcal{I}\in \mathcal{P}(\mathcal{K}), \bigcap_{K\in\mathcal{I}}K\ne\emptyset\Longrightarrow \bigcup_{K\in\mathcal{I}}K\in\mathcal{K}.
\end{displaymath}
\noindent The set $X$ is called the \emph{carrier} of the space $(X,\mathcal{K})$, the set $\mathcal{K}$ is its \emph{connectivity structure}. The elements of  $\mathcal{K}$ are called the \emph{connected subsets} of the space.
The morphisms between two connectivity spaces are the functions which transform connected subsets into connected subsets. They are called the \emph{connectivity morphisms}, or the \emph{connecting maps}\footnote{Though \emph{non-disconnecting maps} would be more accurate.}.
A connectivity space is called \emph{integral} if every singleton subset is connected. The connected subsets with cardinal greater than one will be called the non-trivial connected subsets. A connectivity space is called \emph{finite} if its carrier is a finite set.
\end{df}

 If $X$ is a connectivity space, $\vert X\vert$ will denote its carrier, and $\kappa(X)$ its connectivity structure, so $X=(\vert X\vert,\kappa(X))$.

\begin{rmq} Instead of supposing that the empty set is always a member of connectivity structures, we could suppose without any substantial change  that it is never such a member. But it seems preferable to choose one or the other of those two assumptions, to avoid ``doubling''  the involved categories. 
\end{rmq}

\begin{rmq} Each point of an integral connectivity space belongs to a maximal connected subset. Those subsets are the connected components of the space; they constitute a partition of it.
\end{rmq}

In \cite{Borger:1983}, Börger notes $Zus$ the category of integral connectivity spaces, because of the German word \textit{Zusammenhangsräume}. We propose here to use rather $\mathbf{Cnc}$ to denote the category of connectivity spaces, $\mathbf{Cnct}$ to denote the category of integral connectivity spaces and $\mathbf{fCnct}$ to denote the category of finite integral connectivity spaces.
 
\begin{example}Let $U_T:\mathbf{Top}\to\mathbf{Cnct}$ be the functor whose value is defined on each topological space $(X,\tau)$ as the connectivity space $(X,\mathcal{K})$ with $\mathcal{K}$ the set of connected subsets (in ordinary topological sense) of $(X,\tau)$. Then $U_T$ is not full and  not surjective (up to isomorphism) on objects ; it is faithful but is neither strictly injective nor injective up to isomorphism on objects : for example, if $X=\{a,b\}$, $\tau_1=\{\emptyset, \{a\}, X\}$ and $\tau_2=\{\emptyset, X\}$, then $(X,\tau_1)$ and $(X,\tau_2)$ are not isomorphic but $U_T(X,\tau_1)=U_T(X,\tau_2)$.
\end{example} 

\begin{example}Let $\mathbf{Grf}$ be the topological construct\footnote{\label{construct} Following \cite{joycats:2005}, §5.1, p.\,61, a category of structured sets and structure preserving functions between them is called a \textit{construct}. More precisely, a construct is a concrete category over the category $\mathbf{Set}$ of sets, that is a pair $(\mathbf{A},U)$ where $\mathbf{A}$ is a category and $U:\mathbf{A}\rightarrow \mathbf{Set}$ is a faithful functor (forgetful functor). A \textit{topological} construct is then a construct $(\mathbf{A},U)$ such that the functor $U$ is topological, \textit{i.e.} such that every $U$-structured source $(f_i:E\rightarrow UA_i)_I$ has a unique $U$-initial lift $(\bar{f_i}:A\rightarrow A_i)_I$(see \cite{joycats:2005}, 10.57, p.\,182 and §21.1, p.\,359, and \textit{infra}, the section \ref{Topological categories} of the present article).}
 whose objects are the simple undirected graphs and whose morphisms are the functions which transform edges in edges or in singletons. More precisely, such a graph can be defined as a pair $(X,\mathcal{G})$ with $\mathcal{G}\in\mathcal{Q}_X$ such that 
\[ \{A\in\mathcal{P}_X, \mathrm{card} A=1\}\subseteq \mathcal{G} \subseteq \{A\in\mathcal{P}_X, \mathrm{card} A=2\}, \]
and morphisms $f:(X,\mathcal{G})\to(Y,\mathcal{H})$ are functions $f:X\to Y$ such that $\forall A\in\mathcal{G},f(A)\in \mathcal{H}$. A subset $K$ of such a graph $(X,\mathcal{G})$ is said to be connected if for every pair $(x,x')$ of elements of $K$, there exists a finite path $x=x_0,x_1, \cdots, x_n=x'$ such that each $x_i$ is in $K$ and each $\{x_i,x_{i+1}\}$ is in $\mathcal{G}$. The forgetful functor $U_G:\mathbf{Grf}\to\mathbf{Cnct}$, whose value is defined for each simple undirected graph $(X,\mathcal{G})$ as $(X,\mathcal{K})$ with $\mathcal{K}$ the set of connected subsets of $X$, is a full embedding.\end{example} 

\begin{example}
With each tame link\footnote{A link is called \textit{tame} if it is not \textit{wild}, that is if it is (ambient) isotopic to a polygonal link (or to a smooth link, see \cite{Fox_Crowell:1963}).} $L$ in $\mathbf{R}^3$ or $\mathbf{S}^3$, we associate an integral connectivity space $S_L$ taking the components of the link $L$ as points of $S_L$, the connected subsets of it being defined by the nonsplittable sublinks of $L$. The connectivity structure $\kappa(S_L)$ will be called the \emph{splittability structure} of $L$.
\end{example}

\begin{example}\label{definition of the brunnian spaces} The simplest integral connectivity space which is neither in $U_T(\mathbf{Top})$  nor in $U_G(\mathbf{Grf})$ is the \textit{Borromean} space $\mathbf{B}_3$, defined by $\vert \mathbf{B}_3\vert=3=\{0,1,2\}$ and $\kappa(\mathbf{B}_3)=\mathcal{B}_3$ such that $\mathcal{B}_3^\bullet=\{\vert \mathbf{B}_3\vert\}$. More generally, for each integer $n\in \mathbf{N}$, the $n$-points \textit{Brunnian} space $\mathbf{B}_n$ is the integral connectivity space defined by $\vert\mathbf{B}_n\vert=n$ and $\kappa(\mathbf{B}_n)=\mathcal{B}_n$ such that $\mathcal{B}_n^\bullet=\{\vert \mathbf{B}_n\vert\}$.
The names \textit{Borromean} and \textit{Brunnian} are justified by the fact that the corresponding spaces are the ones associated with the links with same names.
\end{example}

\begin{example} More generally, for each set $X$ and each cardinal $\nu$, there is a unique integral connectivity space whose non-trivial connected subsets are those with cardinal greater than $\nu$.
\end{example}

\begin{example} Let $p$ be an integer. The \emph{hyperbrunnian space} $\mathbf{HB}_p$ is the integral connectivity space such that $\vert\mathbf{HB}_p\vert=\{0,1,\cdots,p-1\}^\mathbf{N}$ and with non-trivial connected subsets all the $K\subseteq{\vert\mathbf{HB}_p\vert}$ for which there exist $k\in\mathbf{N}$ and $a\in\vert\mathbf{HB}_p\vert$ such that $K$ be of the form
\[K=\{x\in\vert\mathbf{HB}_p\vert,\forall n< k, x_n=a_n\}.\]
The space $\mathbf{HB}_3$ will be called the \textit{hyperborromean} space. For each $k\in \mathbf{N}$, the function $\phi_k:\mathbf{HB}_p\to\mathbf{B}_p$ defined by $f(x)=x_k$ is a connectivity morphism. If $p\geq 2$, the function $f:\mathbf{HB}_p\to\mathbf{I}$ defined by \[f(x)=\sum_{n=0}^{n=\infty}\frac{x_n}{p^{n+1}}\] is a surjective connectivity morphism onto $\mathbf{I}=[0,1]$, the connectivity space associated with the usual topological interval $[0,1]$. 
\end{example}

\begin{example} More generally, if $X$ is a set and $(T,\leq)$ is a totally ordered set, we define the integral connectivity space $\mathbf{B}_T(X)$ by $\vert\mathbf{B}_T(X)\vert=X^T$ and $\kappa(\mathbf{B}_T(X))^\bullet=\{K_{f,t}, (f,t)\in X^T\times T\}$ where $K_{f,t}=\{g\in X^T,\forall s\in T, s<t\Rightarrow g(s)=f(s)\}$. Then $\mathbf{B}_p=\mathbf{B}_{\{*\}}(p)$, and $\mathbf{HB}_p=\mathbf{B}_\mathbf{N}(p)$. If $\mathrm{card}(X)\geq 2$, then $\mathbf{B}_T(X)$ is a connected space iff $T$ has a least element. \end{example}

\begin{example} Let $(X,\leq)$ be a totally ordered set. The set of all intervals (of any form) of $X$ constitutes an integral connectivity structure on $X$, called the \textit{order connectivity structure}. In particular, ordinal numbers define connectivity spaces, called the \textit{ordinal connectivity spaces}.
\end{example}

\section{Generation of Connectivity Structures}

\subsection{The Theorem of Generation}

\begin{prop} Let $X$ be a set, and $Cnc_X$  (resp. $Cnct_X$) the set of  connectivity structures on $X$ (resp. the set of integral connectivity structures on $X$). For the order defined by
\[ \mathcal{X}_1\leq \mathcal{X}_2  \Leftrightarrow \mathcal{X}_1\subseteq \mathcal{X}_2, \]
$(Cnc_X,\leq)$ and $(Cnct_X,\leq)$ are complete lattices.
\end{prop}

\noindent \textit{Proof}. These ordered sets have $\mathcal{P}_X$ as a maximal element, and for each nonempty family $(\mathcal{X}_i)_{i\in I}$ of (integral) connectivity structures on $X$, $\bigcap_i \mathcal{X}_i$ is again an (integral) connectivity structure on $X$. 
\begin{flushright}$\square$\end{flushright} 

If $\mathcal{X}_1\leq \mathcal{X}_2$, we say that $\mathcal{X}_1$ is \textit{finer} than $\mathcal{X}_2$, or that $\mathcal{X}_2$ is \textit{coarser} than $\mathcal{X}_1$. $\mathcal{P}_X$, the coarsest structure on $X$, is called the \textit{indiscrete} structure on $X$. The finest connectivity structure contains only the empty set; it is called the \textit{discrete} connectivity structure. The finest integral connectivity structure contains only the empty set and the singletons; it is called the \textit{discrete} integral connectivity structure, or simply the discrete structure.

\begin{rmq} The lattices $Cnc_X$ and $Cnct_X$ are not distributive, unless $X$ has no more than two points. For example, if $X=\{1,2,3\}$ and, for each $i\in X$, $\mathcal{X}_i$ is the integral connectivity structure on $X$ with $(X\setminus \{i\})$ as the only non trivial connected set, then $\bigvee_i (\mathcal{X}_i)=\mathcal{P}_X$, so $\mathcal{B}_3\wedge (\bigvee_i (\mathcal{X}_i))=\mathcal{B}_3$, while $\bigvee_i(\mathcal{B}_3\wedge\mathcal{X}_i)$ is the discrete integral connectivity structure on $X$. \end{rmq}

\begin{df} Let $X$ be a set, and $\mathcal{A}\in \mathcal{Q}_X$ a set of subsets of $X$. The finest  connectivity structure (resp. integral connectivity structure) on $X$ which contains $\mathcal{A}$ is called the \emph{connectivity structure} (resp. \emph{integral connectivity structure}) \emph{generated} by $\mathcal{A}$ and is denoted by $[\mathcal{A}]_0$ (resp. $[\mathcal{A}]$).
\end{df}

Thus, $[\mathcal{A}]_0=\bigwedge\{\mathcal{X}\in Cnc_X, \mathcal{A}\subseteq \mathcal{X} \}$ and $[\mathcal{A}]=\bigwedge\{\mathcal{X}\in Cnct_X, \mathcal{A}\subseteq \mathcal{X} \}$.

\begin{prop}\label{connective function} Let $X$ be a set, $\mathcal{A}$ a set of subsets of $X$, $(Y,\mathcal{Y})$ a connectivity space (resp. integral connectivity space) and $f:X\rightarrow Y$ a function. Then $f$ is a connectivity morphism from $(X,[\mathcal{A}]_0)$ (resp. $(X,[\mathcal{A}])$) to $(Y,\mathcal{Y})$ if and only if $f(A)\in\mathcal{Y}$ for all $A\in \mathcal{A}$.
\end{prop}

\noindent \textit{Proof}. $\{A\in \mathcal{P}_X,f(A)\in \mathcal{Y}\}$ is a connectivity structure on $X$ containing $\mathcal{A}$ and then containing $[\mathcal{A}]_0$ (resp. $[\mathcal{A}]$).\begin{flushright}$\square$\end{flushright} 
\pagebreak[3]

The expression ``\textit{generated structure}''  is justified by the next theorem, in which $\omega_0$ denotes the smallest infinite ordinal.

\begin{thm}[Generation of connectivity structures] \label{Generation of connectivity structures}  Let $X$ be a set and $\mathcal{A}\in \mathcal{Q}_X$ a set of subsets of $X$. Then there exists an ordinal $\alpha_0\leq \omega_0+1$ such that
\[ [\mathcal{A}]_0= \Phi^\alpha(\mathcal{A}) \textrm{ for all }  \alpha\geq\alpha_0, \]
where the $\Phi^\alpha$ are the operators $\mathcal{Q}_X\rightarrow\mathcal{Q}_X$ defined by induction for every ordinal $\alpha$ by
\begin{itemize}
\item $\Phi^0=id_{\mathcal{Q}_X}$,
\item if there is an ordinal $\beta$ such that $\alpha=\beta+1$, then $\Phi^\alpha=\Phi\circ\Phi^\beta$
\item otherwise, for all $\mathcal{U}\in\mathcal{Q}_X$, $\Phi^\alpha(\mathcal{U})=\bigcup_{\beta<\alpha}\Phi^\beta(\mathcal{U})$,
\end{itemize}
and with $\Phi$ the operator  defined for all $\mathcal{U}\in\mathcal{Q}_X$ by
\[ \Phi(\mathcal{U})=\{\emptyset\}\cup\{\bigcup_{A\in \mathcal{E}}A, \mathcal{E}\in \mathcal{L}_\mathcal{U}\},\] 
where $\mathcal{L}_\mathcal{U}=\{\mathcal{E}\in \mathcal{P(U)},\bigcap_{A\in \mathcal{E}}A \neq \emptyset \}.$ 

The integral connectivity structure $[\mathcal{A}]$ generated by $\mathcal{A}$ is  obtained by the same way, adding the singletons of $X$ at any stage of the process.
\end{thm}

\noindent \textit{Proof}. We only have to prove the part of the theorem concerning the generation of connectivity structures, the last claim about \textit{integral} connectivity structure being then obvious.

For every $\mathcal{U}$ and $\mathcal{V}$ in $\mathcal{Q}_X$, we have the three following properties, easy to check:
\begin{itemize}
\item $\mathcal{U}\subseteq\Phi(\mathcal{U})$,
\item $\mathcal{U}\subseteq\mathcal{V}\Rightarrow \Phi(\mathcal{U})\subseteq\Phi(\mathcal{V})$,
\item $\mathcal{U}\in Cnc_X\Leftrightarrow \Phi(\mathcal{U})=\mathcal{U}$.
\end{itemize} 
The first two properties then imply by induction that for every ordinal numbers $\alpha$ and $\beta$ with $\alpha\leq \beta$, one has $\Phi^\alpha(\mathcal{U})\subseteq \Phi^\beta(\mathcal{U})$, and the last two properties imply $\Phi(\mathcal{U})\subseteq [\mathcal{U}]_0$ and, by induction, $\Phi^\alpha(\mathcal{U})\subseteq [\mathcal{U}]_0$ for all ordinal numbers $\alpha$. Then, if for an ordinal number $\alpha_0$ the set $\Phi^{\alpha_0}(\mathcal{A})$ is a connectivity structure on $X$, it coincides with $[\mathcal{A}]_0$. So, to complete the proof, it suffices to verify that the set $\mathcal{C}=\Phi^{\omega_0+1}(\mathcal{A})$ is such a structure, \textit{i.e.} $\Phi(\mathcal{C})=\mathcal{C}$. 
For this, let $\mathcal{W}$ be the set $\Phi^{\omega_0}(\mathcal{A})$, so that $\mathcal{C}=\Phi(\mathcal{W})$. Then $\mathcal{W}$ is stable by union of \textit{finite} famillies with nonempty intersections since $\Phi^{\omega_0}(\mathcal{A})=\bigcup_{n\in \mathbf{N}}\Phi^n(\mathcal{A})$ so every such family is included in $\Phi^n(\mathcal{A})$ for some integer $n$, and its union is again in $\mathcal{W}$. 
Now, let $(S_u)_{u\in U}$ be any family of subsets of $X$ belonging to $\mathcal{C}$ and such that $\bigcap_{u\in U} S_u\neq \emptyset$. We want to verify that $\bigcup_{u\in U} S_u\in\mathcal{C}$.  
For each $u\in U$, $S_u\in \mathcal{C}$ implies that there exists a family $(S_{u,i})_{i\in I_u}$ of subsets of $X$ belonging to $\mathcal{W}$ such that $\bigcap_{i\in I_u}S_{u,i}\neq\emptyset$ and $\bigcup_{i\in I_u}S_{u,i}=S_u$.  
 Let $x$ be an element of $\bigcap_{u\in U} S_u$. For each $u\in U$, there exists an index $i_u \in I_u$ such that $x\in S_{u,i_u}$. For all $u\in U$ and $i\in I_u$, let $T_{u,i}$ be the set $S_{u,i}\cup S_{u,i_u}$. We have $S_{u,i}\in\mathcal{W}$, $S_{u,i_u}\in\mathcal{W}$ and $S_{u,i}\cap S_{u,i_u}\neq\emptyset$ (since $\bigcap_{i\in I_u}S_{u,i}\neq\emptyset$) so $T_{u,i}\in \mathcal{W}$ by the property of $\mathcal{W}$ we emphasized. 
Then $\bigcap_{u\in U, i\in I_U}T_{u,i}\neq\emptyset$, so $\bigcup_{u\in U,i\in I_U} T_{u,i}\in \Phi(\mathcal{W})$, that is $\bigcup_{u\in U} S_u\in \mathcal{C}$.
\begin{flushright}$\square$\end{flushright} 
\pagebreak[3]

\begin{rmq} In the proof behind, the existence of the famillies $(I_u)_{u\in U}$, \linebreak  $((S_{u,i})_{i\in I_u})_{u\in U}$ and $(i_u)_{u\in U}$ depends on the axiom of choice.
\end{rmq} 

\begin{example} Let $X$ be the connectivity space such that $\vert X\vert=\mathbf{R}^2\simeq\mathbf{C}$ and $\kappa(X)=[\mathcal{D}]_0$, where $\mathcal{D}$ is the set of open disks of the Euclidean plane $\mathbf{R}^2$. For $k\in\{1,2,3\}$, let $r_k=(-\frac{1}{2}+\frac{\sqrt{3}}{2}i)^k$ be the cubic roots of unity. For each $(x_0,y_0)=z_0\in \mathbf{C}$, let $(z_n)$ be the sequence of complex numbers defined by the Newton's method for the equation $z^3-1=0$ and with first term $z_0$. If the sequence $(z_n)_{n\in \mathbf{N}}$ converges to $r_k$, we put $f(z_0)=k$, otherwise --- in particular if the sequence $(z_n)$ is defined only for a finite number of terms --- we put $f(z_0)=0$. Then the function $f:X\to \mathbf{B}_4$ defined by this way is a connectivity epimorphism.  Indeed, the three basins of attraction $W_k=f^{-1}(k)$, $k\in\{1,2,3\}$, have the Wada property : their common boundary is the Julia set $W_0=f^{-1}(0)$ (see \cite{Peitgen_Haeseler:1988b}). If $K$ is a nonempty element of the connectivity structure $[\mathcal{D}]_0$, it is open and connected for the usual topology of the plan and then either $K\subset W_k$ for a $k\in\{1,2,3\}$ and $f(K)=\{k\}\in\kappa(\mathbf{B}_4)$, or $K$ intersects $W_{0}$ and then $f(K)=\vert\mathbf{B}_4\vert$ which is again in $\kappa(\mathbf{B}_4)$. 
Note that if we replace $\mathbf{B}_4$ by $(\vert\mathbf{B}_4\vert,\kappa(\mathbf{B}_4)\setminus\{\{0\}\})$, the function $f$ is still a connectivity morphism. Moreover, it is easy to use this function $f$ to define other surjective connectivity morphisms from the same connectivity plane $X$ to the borromean space $\mathbf{B}_3$.
\end{example}

\begin{example} There are several general ways to associate a connectivity space with each (partially) ordered set. We can for example define closed intervals of such a set exactly like in the totally ordered case, and then associate with each ordered set $(S,\leq)$ the connectivity space $(S,[\mathcal{J}]$ with $\mathcal{J}$ the set of closed intervals of $S$. In particular, for each topological construct and each set $X$, we obtain a connectivity space whose points are the structures on $X$.\end{example}

\subsection{Irreducibility}

\begin{df} Let $X$ be a connectivity space. A connected subset $K$ of $\vert X\vert$ is called \emph{reducible} if 
it belongs to the connectivity structure generated by the others, that is 
\begin{displaymath}
K \in [\kappa(X)\setminus \{K\}]_0.
\end{displaymath}
A nonempty connected subset of $\vert X\vert$ is said to be \emph{irreducible} if it is not reducible. The space $X$ is said to be \emph{irreducible} if $\vert X \vert$ is an irreducible connected subset of itself. It is said to be \emph{distinguished} if each of its nonempty connected subsets is irreducible.
\end{df}

\begin{rmq} With the notation of the theorem 
\ref{Generation of connectivity structures}
we have either $\Phi(\kappa(X)\setminus \{K\})=\kappa(X)\setminus \{K\}$, 
 and then $K$ is irreducible, or $\Phi(\kappa(X)\setminus \{K\})=\kappa(X)$. 
 In any case, $[\kappa(X)\setminus \{K\}]_0=\Phi(\kappa(X)\setminus \{K\})$, and $K$ is reducible iff there is a family $\mathcal{E}$ of proper connected subsets $A\subsetneqq K$ such that $\bigcap_{A\in \mathcal{E}}A \neq \emptyset$
and  $K=\bigcup_{A\in \mathcal{E}}A$. 
\end{rmq}

\begin{rmq} A connected singleton is necessarily irreducible.
\end{rmq}

\begin{example} \label{finite irreducible connected subset}
If $X$ is a finite connectivity space, a subset $K$ of $\vert X \vert$ is reducible iff there are two connected subsets $A\subsetneqq \kappa(X)$ and
$B\subsetneqq \kappa(X)$ such that
\begin{displaymath}
K=A\cup B \textrm{ and }
A\cap B \neq \emptyset.
\end{displaymath}
\end{example}

\begin{example} The only irreducible connected subsets of $\mathbf{R}$ are the trivial ones. \end{example}

\begin{example} Brunnian spaces and hyperbrunnian spaces are connected and distinguished spaces. Nevertheless, note that $\mathbf{B}_T(X)$ is not a distinguished space for every set $X$ and every totally ordered set $T$. For example, $\mathbf{B}_{[0,1]}(\{a,b\})$ is not a distinguished space, since $\{f\in \{a,b\}^{[0,1]}, \exists \epsilon \in ]0,1], t<\epsilon \Rightarrow f(t)=a\}$ is a connected subset which is reducible. \end{example}

\begin{df} Let $X$ be a connectivity space. Its \emph{Brunnian closure} is $\overline{X}=(\vert X\vert, \kappa(X)\cup\{\vert X\vert\})$. 
\end{df}

\begin{example} $\mathbf{B}_n$ is the Brunnian closure of the $n$-points discrete integral space. $\mathbf{HB}_n$ is the Brunnian closure of the disjoint union (\textit{cf}. \textit{infra}, section \ref{disjoint union of spaces}) of $n$ copies of itself.\end{example}

The next proposition is obvious.

\begin{prop} If $X$ is a nonempty irreducible space, then $(\vert X\vert, \kappa(X)\setminus \{X\})$ is a connectivity space. If $X$ is a non-connected connectivity space, then $\overline{X}$ is an irreducible connected space. 
\end{prop}

Because of the next proposition, the notion of irreducibility will play a fundamental role in the case of finite connectivity spaces.

\begin{prop}\label{characterization of finite spaces by irreducible subsets} A connectivity structure on a given finite set is characterised by the set of the irreducible connected subsets, which is the minimal set of subsets which generates this structure.
\end{prop}

\noindent\textit{Proof}. For any connectivity space $X$, let $\iota(X)$ denote the set of the irreducible connected subsets of $X$. Then, for any $\mathcal{A}\in \mathcal{Q}_X$ such that $[\mathcal{A}]_0=\kappa(X)$, one has $\mathcal{A}\supseteq \iota(X)$ since, by construction, each set $C\in [\mathcal{A}]_0$ which is not in $\mathcal{A}$ is reducible. On the other hand, an easy induction shows that, for every integer $k$, every reducible connected subset of $X$ with cardinal smaller than $k$ is an element of $[\iota(X)]_0$. Thus, if $X$ is finite, $\kappa(X)=[\iota(X)]_0$.
\begin{flushright}$\square$\end{flushright}

\subsection{Connectivity Spaces and Hypergraphs}

A hypergraph is a set of vertices endowed with a set of nonempty sets of vertices, these sets of vertices being considered as generalized edges, the so-called \textit{hyperedges}. There is some similarity between hypergraphs and connectivity spaces --- for example it is possible to consider Borromean structures in both cases --- but 

\begin{itemize}
\item the union of two hyperedges with a nonempty intersection is not necessarily an hyperedge, so hyperedges are not the same as connected subsets,
\item the union of two hyperedges with a nonempty intersection can be an hyperedge, so hyperedges are not the same as irreducible connected subsets.
\end{itemize}

To clarify the relation between the two concepts, let us consider the category $\mathbf{HypG}$ of hypergraphs, that is the category 
whose objects are the pairs $(X,\mathcal{H})$ with $X$ a set and $\mathcal{H}\in\mathcal{Q}_X$ a set whose elements are called hyperedges, and whose morphisms $f:(X,\mathcal{H})\to(X',\mathcal{H'})$ are functions $X\to X'$ which preserve hyperedges : $H\in\mathcal{H}\Rightarrow f(H)\in\mathcal{H'}$. Then the proposition \ref{connective function} implies 

\begin{cor} \label{relation between connectivity and hypergraphs}
The category $\mathbf{Cnc}$ is concrete on $\mathbf{HypG}$ with a forgetful functor admiting as a left adjoint the functor $\mathbf{HypG}\to\mathbf{Cnc}$ which associates with each hypergraph $(X,\mathcal{H})$ the space whose connectivity structure is generated by $\mathcal{H}$, \textit{i.e.} $(X,[\mathcal{H}]_0)$, and with each morphism itself as a connectivity morphism. Similarly, the generation of integral connectivity structures $[\mathcal{H}]$ from sets $\mathcal{H}\in \mathcal{Q}_X$ defines a left adjoint to the forgetful functor $\mathbf{Cnct}\to\mathbf{HypG}$, and the situation is the same between finite hypergraphs and finite connectivity spaces.
\end{cor}

\section{Limits and Colimits}

\subsection{Categories with Lattices of Structures}
\label{Topological categories}

Let $\mathbf{JCPos}$ be the category of complete (small) lattices and join-preserving maps. If $S$ is a functor from a category $\mathbf{X}$ to $\mathbf{JCPos}$, $S(X)$  or $S_X$ will denote the lattice associated by $S$ with an object $X$, and (while it is unambigous)
$f_!$ the map between lattices associated by $S$ with a morphism $f$. The elements of the lattice $S_X$ will be called the $S$-structures on $X$.

\begin{df} If $\mathbf{X}$ is a category and $S:\mathbf{X}\to\mathbf{JCPos}$ a functor, the \emph{category with lattices of structures associated with $S$}, which will equally be called the
\emph{category structured by $S$}, is the category we denote $\mathbf{X}_S$, whose objects are the pairs $(X,\mathcal{X})$ with $X$ an object of $\mathbf{X}$ and $\mathcal{X}\in S_X$ a $S$-structure, and whose morphisms $f:(X,\mathcal{X})\to(Y,\mathcal{Y})$ are $\mathbf{X}$-morphisms $X\to Y$ such that $f_!(\mathcal{X})\leq \mathcal{Y}$ in the lattice $S_Y$.
\end{df}

In the category $\mathbf{X}_S$, spaces $(X,1_{S_X})$ are called \textit{indiscrete} spaces, and spaces $(X,0_{S_X})$ are called \textit{discrete} spaces. If, in the lattice $S_X$, we have $\mathcal{X}\leq \mathcal{X}'$, then the structure $\mathcal{X}$ is said to be \textit{finer} than $\mathcal{X}'$ and the latter is said to be \textit{coarser} than the former.


\begin{rmq} 
An equivalent definition is given by considering \emph{contravariant} functors from the basis category $\mathbf{X}$ to the category $\mathbf{MCPos}$ of complete (small) lattices and meet-preserving maps: an object of the category defined by such a functor $T$ is a pair $(X,\mathcal{X})$ with $\mathcal{X}\in T_X$, and a morphism $f:(X,\mathcal{X})\to(Y,\mathcal{Y})$ is a $\mathbf{X}$-morphism $f:X\to Y$ such that $\mathcal{X}\leq f^*(\mathcal{Y})$, where $f^*=T(f)$. Then, for each covariant $S:\mathbf{X}\to\mathbf{JCPos}$, there is a \emph{contravariant associated functor} $T$ defining by this way the  category we noticed $\mathbf{X}_S$. This functor  $T$ is defined on objects $X$ by $T_X=S_X$ and on $\mathbf{X}$-morphisms $f:X\to Y$ by $T(f)=f^*$ with, for each $\mathcal{Y}\in T_Y$, 
\[  f^*\mathcal{Y}=\bigvee \{\mathcal{X}\in T_X,f_!\mathcal{X}\leq\mathcal{Y}\}.\]
\end{rmq}

In the next proposition, we use the definition of a topological category given in \cite{joycats:2005} : a topological category on $\mathbf{X}$ is a concrete category $U:\mathbf{A}\rightarrow\mathbf{X}$ (that is, a faithful functor $U$), such that every $U$-source $(X\rightarrow UA_i)_{i\in I}$ in $\mathbf{X}$ has a unique $U$-initial lift $(A\rightarrow A_i)_{i\in I}$ in $\mathbf{A}$. 

\begin{prop} A category is a small-fibred topological one if and only if it is a category with lattices of structures. More precisely : 
\begin{itemize}
\item For each functor $S:\mathbf{X}\rightarrow\mathbf{JCPos}$, the functor $U:\mathbf{X}_S\rightarrow \mathbf{X}$ defined by $U(X,\mathcal{X})=X$ and $Uf=f$ is a small-fibred topological category.
\item Each small-fibred topological category $U:\mathbf{A}\rightarrow\mathbf{X}$ is isomorphic to the category $\mathbf{X}_S$ with $S$ the functor defined for each object $X$ of $\mathbf{X}$ by the fibre $S_X=\{A\in \mathbf{A}, UA=X\}$ with the usual order (i.e. $A_1\leq A_2$ iff $id_X$ has a lift $A_1\rightarrow A_2$), and for each arrow $f:X\rightarrow Y$ in $\mathbf{X}$ and each $A\in S_X$ by $f_{!}(A)=\wedge \{B\in S_Y, f \textrm{ has lift } A\rightarrow B\}$.
\end{itemize}
\end{prop}

\noindent \textit{Proof}. Let $S:\mathbf{X}\rightarrow\mathbf{JCPos}$ be any functor. The functor $U:\mathbf{X}_S\rightarrow \mathbf{X}$ defined by $U(f:(X,\mathcal{X})\rightarrow(Y,\mathcal{Y}))=(f:X\rightarrow Y)$ 
is trivially faithful, its fibres are the sets $S_X$, and it is topological : each 
$U$-source $(f_i:X\rightarrow UA_i)_{i\in I}$
 has a unique $U$-initial lift, that is 
$(f_i:(X,\mathcal{X}_0)\rightarrow A_i)_{i\in I}$,
 where $\mathcal{X}_0$ is the coarsest $S$-structure on $X$ such that all $f_i$ be (have lifts as) $\mathbf{X}_S$-morphisms, that is
 $\mathcal{X}_0=\wedge_i{f^*_i(\mathcal{Y} _i)}$
 where $\mathcal{Y} _i$ is the $S$-structure of $A_i$ and, for each $f:X\rightarrow Y$
 and each 
$\mathcal{Y}\in S_Y$,$f^*(\mathcal{Y})$
 is the coarsest $S$-structure $\mathcal{X}$ on $X$ such that $f$ be a $\mathbf{X}_S$-morphism 
$(X,\mathcal{X})\rightarrow(Y,\mathcal{Y})$,
 that is $f^*(\mathcal{Y})=\vee\{\mathcal{X}\in S_X, f_!(\mathcal{X}) \leq \mathcal{Y}\}$.

On the other hand, let now $U:\mathbf{A}\rightarrow\mathbf{X}$ be a topological category with small fibres. One knows (see \cite{joycats:2005}) that such fibres $S_X$ are then complete lattices. We can remark also that, for a given $f:X\rightarrow Y$ in $\mathbf{X}$ and an object $A\in S_X$, the set $\{B\in S_Y, f \textrm{ has a lift } A\rightarrow B\}$ is nonempty, because $Y$ has an indiscrete lift. Then $f_!$ is well-defined as a function. Now, if $(A_i)_{i\in I}$ is any family in the fibre $S_X$, and $B\in S_Y$ is such that $f:X\rightarrow Y$ has a lift $\vee_i A_i\rightarrow B$, then $id_X$ has a lift $A_i\rightarrow\vee_i A_i$ for each $i$, so $f$ has a lift $A_i\rightarrow B$ for each $i$. On the other hand, if $f:X\rightarrow UB$ has a lift $A_i\rightarrow B$ for each $i$, then $\forall i\in I, A_i\leq A$, where the $U$-initial lift of $f$ is $A\rightarrow B$; but $\vee_i A_i\leq A$, so $id_X$ has a lift $\vee_i A_i\rightarrow A$ and $f$ has a lift $\vee_i A_i\rightarrow B$. Thus, for a given $f:X\rightarrow Y$ and a given family $(A_i)_{i \in I}$ in $S_X$, we have 
$\{B\in S_Y, f \textrm{ has a lift }\vee_i A_i\rightarrow B\}=\{ B\in S_Y, \forall i\in I, f \textrm{ has a lift } A_i\rightarrow B\}$.
Let $\beta_i=f_!(A_i)=\wedge\{ B\in S_Y, f \textrm{ has a lift } A_i\rightarrow B\}$. Then
\[ f_!(\vee_i A_i)
=\wedge \{ B\in S_Y, \forall i\in I, f \textrm{ has a lift } A_i\rightarrow B\} \]
\[ = \wedge \{ B\in S_Y, \forall i\in I, B\geq \beta_i\}=\vee_i\beta_i, \]
so $f_!(\vee_i A_i)=\vee_i f_!(A_i)$ : $f_!$ is a $\mathbf{JCPos}$-morphism, and the functor $S$ is well-defined. It is then easy to verify that the functor $\mathbf{A}\rightarrow \mathbf{X}_S$ defined by 
\[ (f:A\rightarrow B)\mapsto (Uf:(UA,A)\rightarrow (UB,B))\]
is an isomorphism of categories, with inverse 
\[ (f: (X,\mathcal{X})\rightarrow (Y,\mathcal{Y})) \mapsto (\tilde{f}:\mathcal{X}\rightarrow\mathcal{Y}),\] where $\tilde{f}$ is the lift of $f$, which exists since $f_!(\mathcal{X})\leq\mathcal{Y}$.
\begin{flushright}$\square$\end{flushright} 
\pagebreak[3]

By the proposition $21.15$, the theorem $21.16$ and the corollary $21.17$ of \cite{joycats:2005}, we have then

\begin{cor} \label{Properties of categories with lattices of structures} If  $\mathbf{X}$ denotes the category $\mathbf{Set}$ of sets (resp. the category $\mathbf{fSet}$ of finite sets), $S:\mathbf{X}\to\mathbf{JCPos}$ any functor, $T:\mathbf{X}^{op}\to\mathbf{MCPos}$ the contravariant functor associated with $S$ and $U:\mathbf{A}=\mathbf{X}_S\to\mathbf{X}$ the construct\footnote{See \textit{supra} the note \ref{construct}.} (resp. ``finitely''  construct) defined by $S$, then the following hold
\begin{enumerate}
\item $\mathbf{A}$ is (co)complete (resp. finitely (co)complete), 
\item $U$ has a left adjoint $O$ (the discrete structure) and a right adjoint $I$ (the indiscrete structure) : $O\dashv U\dashv I$, so $U$ preserves (co)limits,
\item the limit $(l_i:L\to D_i)_{i\in \mathbf{I}}$ of a small (resp. finite) diagram $D:\mathbf{I}\to\mathbf{A}$ is the initial lift of the underlying limit in $\mathbf{X}$, that is: if $(l_i:\vert L\vert\to UD_i)_{i\in \mathbf{I}}$ is the limit of $UD$, then $L=(\vert L\vert,\bigwedge_{i\in\mathbf{I}}l_i^*(\mathcal{X}_i))$, where $D_i=(X_i,\mathcal{X}_i)$ and $l_i^*=T(l_i)$,
\item colimits are given in the same way, as final lifts: if $(c_i:\vert C\vert\leftarrow UD_i)_{i\in \mathbf{I}}$ is the colimit of $UD$, then the colimit of $D$ in $\mathbf{A}$ is $(c_i:C\leftarrow D_i)_{i\in \mathbf{I}}$ with $C=(\vert C\vert,\bigvee_{i\in\mathbf{I}}{c_i}_!(\mathcal{X}_i))$, where $D_i=(X_i,\mathcal{X}_i)$ and ${c_i}_!=S(c_i)$,
\item $\mathbf{A}$ is wellpowered and cowellpowered,
\item $\mathbf{A}$ is an $(Epi, Extremal MonoSource)$-category, 
\item $\mathbf{A}$ has regular factorizations, \textit{i.e.} is an $(RegEpi, MonoSource)$-category (and thus is, in particular, a $(RegEpi,Mono)$-category),
\item in $\mathbf{A}$, the classes of embeddings (i.e. initial monomorphisms), of extremal monomorphisms and of regular monomorphisms coincide,
\item in $\mathbf{A}$, the classes of quotient morphisms (i.e. final epimorphisms), of extremal epimorphisms and of regular epimorphisms coincide,
\item $\mathbf{A}$ has separators and coseparators.
\end{enumerate}
\end{cor}

\begin{example}
Let $\mathcal{P}:\mathbf{Set}\rightarrow \mathbf{JCPos}$ be the (covariant) functor which associates with each set the complete lattice of its subsets. For any functor $T:\mathbf{X}\rightarrow \mathbf{Set}$, the category $\mathbf{X}_{\mathcal{P}T}$ structured by the functor $\mathcal{P}\circ T : \mathbf{X}\rightarrow \mathbf{JCPos}$ coincides with the topological category $\mathbf{Spa}(T)$ of $T$-spaces on $\mathbf{X}$ (\cite{joycats:2005}, p.\,76). Thus, the ``functor-structured categories'' $\mathbf{Spa}(T)$ are special cases of the categories structured by functors $\mathbf{X}\to\mathbf{JCPos}$. In particular, for $T=\mathcal{P}$, we obtain $\mathbf{Spa}(\mathcal{P})=\mathbf{Set}_{\mathcal{Q}}=\mathbf{HypG}$.
\end{example}

\subsection{(Co)limits in the Categories of Connectivity spaces}

In \cite{Borger:1983}, Börger showed that 

\begin{prop} $\mathbf{Cnct}$ is a topological category. It is not cartesian closed.
\end{prop}

It is easy to check that, as a category with lattice of structures, $\mathbf{Cnct}$ is defined by the covariant functor $Cnct:\mathbf{Set}\to \mathbf{JCPos}$ such that $Cnct_X$ is the lattice of all integral connectivity structures on $X$ and, for every $f:X\to X'$, $Cnct(f)=f_!$ is the $\mathbf{JCPos}$-morphism $Cnct_X\to Cnct_{X'}$ such that, for all $\mathcal{K}\in Cnct_X$,
\begin{equation} \label{expression of f!}
f_!(\mathcal{K})=[\{f(K), K\in \mathcal{K}\}].
\end{equation}
Equivalently, the contravariant definition of $\mathbf{Cnct}$ is given, for all $\mathcal{K'}\in Cnct_{X'}$, by
\begin{equation}  \label{expression of f*}
f^*(\mathcal{K'})=\{K\in \mathcal{P}_X, f(K)\in \mathcal{K'}\}.
\end{equation}

The same formulas hold on $\mathbf{fSet}$, defining a functor $fCnct$ such that $\mathbf{fCnct}=\mathbf{fSet}_{fCnct}$, which is thus a topological category on $\mathbf{fSet}$. For $\mathbf{Cnc}$, it suffices to use $[\{f(K), K\in \mathcal{K}\}]_0$ instead of $[\{f(K), K\in \mathcal{K}\}]$ in the expression of $f_!$ to define a functor $Cnc$ such that $\mathbf{Cnc}=\mathbf{Set}_{Cnc}$, which is thus a topological construct\footnote{$\mathbf{Cnc}$ is not \textit{well-fibred}, so it is not a topological category according to the definition given in 1983 by Herrlich \cite{Herrlich:1983}, but, as we said, we use here the less restrictive definition finally retained by Herrlich, Ad\'{a}mek and Strecker in \cite{joycats:2005}.}.

From the formula (\ref{expression of f!}) and the corollary \ref{Properties of categories with lattices of structures}, we deduce that the connectivity structure $\kappa(C)$ of the colimit $C$ of a small diagram $D:\mathbf{I}\to\mathbf{Cnct}$ is given by $\kappa(C)=\bigvee_{i\in \mathbf{I}}[\{c_i(K), K\in\kappa(D_i)\}]$ and then
\begin{equation}\label{expression of colimits} 
\kappa(C)=[\{c_i(K), i\in \mathbf{I}, K\in\kappa(D_i)\}],
\end{equation}
where the $c_i:\vert D_i\vert \to \vert C\vert $ are the coprojections.
The same formula holds for colimits of finite diagrams in $\mathbf{fCnct}$, and, using $[-]_0$ instead of $[-]$, for small diagrams in $\mathbf{Cnc}$. 

From the formula (\ref{expression of f*}), one likewise deduces the connectivity structure $\kappa(L)$ of the limit $L$ of a small diagram $D:\mathbf{I}\to\mathbf{Cnct}$,
\begin{equation}\label{expression of limits} 
\kappa(L)=\bigcap_{i\in \mathbf{I}}\{K\in \mathcal{P}_{\vert L\vert},l_i(K)\in\kappa(D_i)\},
\end{equation}
where the $l_i:\vert L\vert\to\vert D_i\vert$ are the projections.
The same formula holds for limits of small diagrams in $\mathbf{Cnc}$ and of finite diagrams in $\mathbf{fCnct}$. 

For example, the cartesian product $C_1\times C_2$ of two connectivity spaces is characterised by $\vert C_1\times C_2\vert = \vert C_1\vert\times \vert C_2\vert$ and \[\kappa (C_1\times C_2)= \{A\in \mathcal{P}(\vert C_1\vert\times \vert C_2\vert), \pi_i(A)\in \kappa(C_i) \textrm{ for } i\in\{1,2\}\},\] where the $\pi_i$ are the projections, whereas the coproduct, or disjoint union\label{disjoint union of spaces}, satisfies $\vert C_1\amalg C_2\vert = \vert C_1\vert\amalg \vert C_2\vert $ and $\kappa(C_1\amalg C_2) = \kappa( C_1)\amalg \kappa( C_2) $.

With those formulas, it is easy to check that none of the three categories considered here is cartesian closed. It suffices to exhibit a colimit which is not preserved by a product, and this can be done simultaneously in the three categories. For example, let $\{a,*,b\}$ be a set with three distincts elements, $A_u$ be the indiscrete connectivity space defined for each $u\in\{a,b\}$ by its carrier $\vert A_u\vert=\{*,u\}$, and $B$ the space with carrier $\{1,2,3\}$ and with structure $[\{\{1,2\},\{2,3\}\}]$. Then, in each of the categories concerned, the colimit $C$ of the diagram $A_a\hookleftarrow\{*\}\hookrightarrow A_b$ (with arrows the inclusions) is $C=(\{a,*,b\},[\{\{a,*\},\{*,b\}\}])$, its product $C\times B$ with $B$ is the cartesian product $\{a,*,b\}\times\{1,2,3\}$ endowed with the integral connectivity structure including all subsets having their two projections connected. For example, the set $\{(a,1),(*,3),(b,2)\}$ is connected in $C\times B$; but it is easy to verify that the same set is not connected in the colimit of the diagram $A_a\times B \hookleftarrow \{*\}\times B \hookrightarrow A_b\times B$. Thus, in each of the categories considered, the endofunctor $-\times B$ does not preserve colimits. We thus proved 

\begin{prop} $\mathbf{Cnc}$ and $\mathbf{fCnct}$ are topological categories; they are not cartesian closed.
\end{prop}

\subsection{Quotients and Embeddings}

This section gives trivial but useful consequences of the corollary \ref{Properties of categories with lattices of structures} and of the formulas (\ref{expression of f!}) and (\ref{expression of f*}).

\begin{prop} 

In $\mathbf{Cnct}$ and $\mathbf{fCnct}$ (resp. $\mathbf{Cnc}$), a morphism $f:A \to B$ is a \emph{regular epimorphism} iff $\vert f\vert$ is surjective and $\kappa(B)=[f(\kappa(A))]$ (resp. $\kappa(B)=[f(\kappa(A))]_0$). 
In $\mathbf{fCnct}$, $\mathbf{fCnct}$ and $\mathbf{Cnc}$, a morphism $f:A \to B$ is a \emph{regular monomorphism} iff $\vert f\vert$ is injective and $\kappa(A)=\{K\in \mathcal{P}_{\vert A\vert}, f(K)\in \kappa(B)\}$.
\end{prop}

Now, in every topological construct, a regular epimorphism, \textit{i.e.} a coequalizer, is the same as a quotient morphism, \textit{i.e.} a final morphism which is surjective as a function, and can also be viewed as (the unique final lift of) the canonical map associated with an equivalence relation. This remark results in the definition of the quotient of a connectivity space by an equivalence relation.

\begin{df}[Quotient by an equivalence relation] If $C$ is a connectivity space and $\sim$ is an equivalence relation on $\vert C\vert$, the \emph{quotient space} $C/\sim$ is defined by $\vert C/\sim\vert=\vert C\vert/\sim$ and 
\begin{equation}
\kappa(C/\sim)= s_!(\kappa(C))=[s(\kappa(C))]_0
\end{equation}
where $s$ is the canonical map $s:\vert C \vert \twoheadrightarrow \vert C \vert/\sim$. In particular, if $T$ is a subset of $\vert C\vert$, $C/T$ denotes the space $C/\sim_T$.
\end{df}

\begin{rmq} Note that if $C$ is an integral connectivity space, then for any surjective map $s:\vert C \vert \twoheadrightarrow Y$ we have $[s(\kappa(C))]_0=[s(\kappa(C))]$.
\end{rmq}

Likewise, in every topological construct, a regular monomorphism, \textit{i.e.} an equalizer, is the same as an embedding, \textit{i.e.} an initial morphism which is injective as a function, and can also be viewed as (the unique initial lift of) the inclusion map of a subspace. This leads to the definition of the connectivity structure induced by a connectivity space on a subset of its carrier.

\begin{df}[Structure induced on a subset] If $C$ is a connectivity space and $S$ is a subset of $\vert C\vert$, the \emph{connectivity space induced} on $S$ by $C$ is the space $C_{\vert S}$ defined by $\vert C_{\vert S}\vert=S$ and 
\begin{equation}
\kappa(C_{\vert S})= i^*(\kappa(C))=\mathcal{P}_S\cap\kappa(C)
\end{equation}
where $i$ is the inclusion map $i:S \hookrightarrow \vert C \vert$.
\end{df}

\section{Tensor Product of Connectivity Spaces}

The formula (\ref{expression of limits}) suggests that the cartesian product of connectivity spaces is in some way ``too coarse'' to be really useful in algebra. For example, let $\mathbf{N}$ be the set of natural numbers with the integral connectivity structure generated by the subsets $\{n,n+1\}$; it is easy to check that the addition $+:\mathbf{N}^2\to\mathbf{N}$ is not a connectivity morphism (when $\mathbf{N}^2$ is endowed with the cartesian square structure of $\mathbf{N}$). Likewise for the addition of real numbers. This section presents a more interesting connectivity product than the cartesian one for algebraic structures.
 
Let $X_i$ ($i=1,2$) and $Y$ be connectivity spaces. For each $x_1\in \vert X_1\vert$ (resp. $x_2\in \vert X_2\vert$), we denote by $f(x_1,-)$ (resp. $f(-,x_2)$) the partial function associated with a given function $f:\vert X_1\vert\times \vert X_2\vert\to \vert Y\vert$.

\begin{df} A function $f:\vert X_1\vert\times \vert X_2\vert\to \vert Y\vert$ is said to be \emph{partially connecting} from $X_1\times X_2$ to $Y$ if $f(x_1,-):X_2\to Y$ and $f(-,x_2):X_1\to Y$ are connectivity morphisms for all $x_1\in \vert X_1\vert$ and all $x_2\in\vert X_2\vert$.
\end{df}

\begin{df} The \emph{connectivity tensor product} $X_1\boxtimes X_2$ of two connectivity spaces $X_i$ ($i=1,2$) is the space with carrier $\vert X_1\boxtimes X_2 \vert=\vert X_1\vert\times \vert X_2\vert$ and with connectivity structure $\kappa(X_1\boxtimes X_2)=[\{K_1\times K_2, (K_1,K_2)\in \kappa(X_1)\times\kappa(X_2)\}]_0$.
\end{df}

For every connectivity space $X_i$, $\kappa(X_1\boxtimes X_2)$ is a finer connectivity structure on the set $\vert X_1\vert\times \vert X_2\vert$ than the one given by the connectivity cartesian product, since $K_1\times K_2\in \kappa(X_1\times X_2)$ for each connected subsets $K_1$ and $K_2$. Thus, $id:X_1\boxtimes X_2\to X_1\times X_2$ is a bijective connectivity morphism (but it is of course not an isomorphism in general). If $X_1$ and $X_2$ are \textit{integral} connectivity spaces, then its inverse function, that is the function from $X_1\times X_2$ to $X_1\boxtimes X_2$ defined by $\tau(x_1,x_2)=(x_1,x_2)$, is a partially connecting function.

\begin{thm} Let $X_1$ and $X_2$ be integral connectivity spaces, $Y$ a connectivity space, and $f:\vert X_1\vert\times \vert X_2\vert\to \vert Y\vert$ a function. Then $f$ is a partially connecting function  from $X_1\times X_2$ to $Y$ if and only if it is a connectivity morphism from $X_1\boxtimes X_2$ to $Y$, \emph{i.e.} there exists a unique connectivity morphism $\tilde{f}:X_1\boxtimes X_2\to Y$ such that $\tilde{f}\circ \tau=f$.
\end{thm}
\noindent \textit{Proof}. If $\tilde{f}$ is a connectivity morphism, then $\tilde{f}\circ \tau=f$ is a partially connecting function since $\tau$ is such a function. On the other hand, let $f$ be a partially connecting function from $X_1\times X_2$ to $Y$. Unicity of $\tilde{f}$ being obvious, since necessarily $\tilde{f}(x_1,x_2)=f(x_1,x_2)$, it suffices to check that this function is a connectivity morphism on $X_1\boxtimes X_2$. Then, according to the proposition \ref{connective function}, it suffices to check that for every $K_i\in\kappa(X_i)$, $f(K_1\times K_2)\in\kappa(Y)$. Let $K_1\times K_2$ be such nonempty subset of $\vert X_1\vert\times \vert X_2\vert$, and let $x_1^0\in K_1$. $f$ being partially connecting, the sets $V=\{f(x_1^0,x_2), x_2\in K_2\}$ and $H_{x_2}=\{f(x_1,x_2), x_1\in K_1\}$ are, for all $x_2\in K_2$, in  $\kappa(Y)$. So are the sets $V\cup H_{x_2}$ (as $V\cap H_{x_2}\neq\emptyset$), and $\bigcup_{x_2\in K_2}(V\cup H_{x_2})$; that is: $\tilde{f}(K_1\times K_2)\in\kappa(Y)$.\begin{flushright}$\square$\end{flushright} 

\begin{example} Let $f:\mathbf{R}_+^2\to \mathbf{R}$ defined by
\begin{itemize}
\item $f(0,0)=0$,
\item for all $x$ and $y$, $f(x,y)=f(y,x)$,
\item $\forall x>0, \forall y\in [0,x], f(x,y)=y/x$.
\end{itemize}
Then $f$ is a partially connecting map since it is ``partially continuous'', but it is not continuous, and neither $\Delta=\{(x,x),x\geq 0\}$ nor $f(\Delta)=\{0,1\}$ are connected subsets of, respectively, $\mathbf{R}_+ \boxtimes \mathbf{R}_+$ and $\mathbf{R}$.
\end{example}

Note that for each integral connectivity space $X$, one has an endofunctor $X\boxtimes - : \mathbf{Cnct} \to \mathbf{Cnct}$ defined for each integral connectivity space $Y$ by $X\boxtimes Y$ and for each connectivity morphism $g:Y_1\to Y_2$ between integral connectivity spaces by $(X\boxtimes g)(x,y_1)=(x,g(y_1))$. 

Now, let us define another endofunctor on $\mathbf{Cnct}$. For every subset $M$ of the set $Hom(X,Y)$ of connectivity morphisms from a connectivity space $X$ to a connectivity space $Y$, and for every subset $A$ of the set $\vert X \vert$, let $\langle M,A\rangle$ denotes $\bigcup_{f\in M} f(A)$. Then, for each integral connectivity space $X$, there is an endofunctor $\mathbf{Cnct}(X,-):\mathbf{Cnct}\to\mathbf{Cnct}$ defined  for every integral connectivity space $Y$ by
\begin{itemize}
\item $\vert \mathbf{Cnct}(X,Y)\vert=Hom(X,Y)$,
\item $\kappa(\mathbf{Cnct}(X,Y))=\{M\in\mathcal{P}(Hom(X,Y)),\forall K\in\kappa(X),\langle M,K\rangle\in\kappa(Y)\}$,
\end{itemize}
and for every connectivity morphism $g:Y_1\to Y_2$ by $\mathbf{Cnct}(X,g)=g_*$ such that
\[ 
\forall \varphi\in \mathbf{Cnct}(X,Y_1), g_*(\varphi)=g\circ \varphi.
\]

\begin{rmq} A set $M$ of connectivity morphisms between two integral connectivity spaces $X$ and $Y$ is connected, that is belongs to $\kappa(\mathbf{Cnct}(X,Y))$, if (and only if) for all $x\in X$, $\langle M,\lbrace x\rbrace\rangle\in\kappa(Y)$. Indeed, if this condition is satisfied, then for every nonempty connected subset $K$ of $X$ and any $x\in K$, one has $\langle M,K\rangle=\bigcup_{f\in M}(f(K)\cup \langle M,\lbrace x\rbrace\rangle)\in\kappa(Y)$.
\end{rmq}

\begin{thm} For every integral connectivity space $X$, the endofunctor \linebreak $X\boxtimes -$ is left adjoint to the endofunctor $\mathbf{Cnct}(X,-)$. Thus, $(\mathbf{Cnct},\boxtimes)$ is a closed symmetric monoidal category.
\end{thm}

\noindent \textit{Proof}. The product $\boxtimes$ is obviously symmetrical. Let $X$, $Y$ and $Z$ be integral connectivity spaces. For every connectivity morphism $\psi:X\boxtimes Y\to Z$, one has a morphism $\rho(\psi):Y\to \mathbf{Cnct}(X,Z)$ defined for all $y\in Y$ by $\rho(\psi)(y)=\psi(-,y)$. Then $\rho$ is clearly a bijection between the sets $Hom(X\boxtimes Y,Z)$ and $Hom(Y,\mathbf{Cnct}(X,Z))$, and it is natural since for all integral connectivity spaces $Y$, $Y'$, $Z$ and $Z'$ and for all connectivity morphisms $u:Y\to Y'$, $v:Z\to Z'$ and $\psi:X\boxtimes Y'\to Z$, one has $\rho(v\circ \psi\circ (X\boxtimes u))=\rho((x,y)\mapsto v(\psi(x,u(y))))=(y\mapsto v\circ \psi(-,u(y))=\mathbf{Cnct}(X,v)\circ \rho(\psi)\circ u$.
\begin{flushright}$\square$\end{flushright} 

\section{Homotopy}

Let $\overrightarrow{I}$ be a triple $(I,0,1)$ with $I$ a nonempty integral connectivity space, and $0$ and $1$ some elements of $\vert I\vert$. In particular, let $\mathbf{I}$ be the connectivity space associated with the usual topological space $[0,1]$, and $\overrightarrow{\mathbf{I}}=(\mathbf{I},0,1)$. 

\begin{df} [Homotopy] Let $X$ and $Y$ be integral connectivity spaces, and $f, g : X\to Y$ some connectivity morphisms. The function $g$ is said to be $\overrightarrow{I}$-homotopic to $f$ provided there exists a connectivity morphism \[h:I\to \mathbf{Cnct}(X,Y)\] such that $h(0)=f$ and $h(1)=g$. In particular, in the case of $\overrightarrow{I}=\overrightarrow{\mathbf{I}}$, $g$ is simply said to be homotopic to $f$.
\end{df}

We denote by $f\sim g$ the homotopy relation between connectivity morphisms. Like in the topological case, it is obviously an equivalence relation. The adjoint situation $(X\boxtimes -)\dashv\mathbf{Cnct}(X,-)$ leads to an alternative definition of homotopy for connectivity morphisms.

\begin{df} [Alternative definition of homotopy] Let $X$ and $Y$ be integral connectivity spaces. A function $g:X\to Y$ is $\overrightarrow{I}$-homotopic to $f:X\to Y$ provided there exists a connectivity morphism $h:I\boxtimes X\to Y$ such that $h(0,-)=f$ and $h(1,-)=g$, that is a function $h:I\times X\to Y$ such that
\begin{itemize}
\item $h(0,-)=f$ and $h(1,-)=g$,
\item $\forall t\in I, \forall K\in \kappa(X), h(t,K)\in\kappa(Y)$,
\item $\forall D\in \kappa(I), \forall x\in X, h(D,x)\in\kappa(Y)$.
\end{itemize}
\end{df}

\begin{df} [Contractibility] An integral connectivity space $X$ is said to be \emph{contractible} provided the identity map $id:X\to X$ of the space be homotopic to a constant map $c:X\to X$.
\end{df}

\noindent\textit{Examples}. The connectivity space associated with the usual topological circle $S^1=\{e^{i\theta},\theta\in[0,2\pi]\}\subset \mathbf{C}$ is contractible. Indeed, the function $h:\mathbf{I}\times S^1\to S^1$ defined by 
\begin{itemize}
\item for $t\in [0,1[$ and $z\in S^1$, $h(t,z)=z.e^{i\frac{t}{1-t}}$,
\item $\forall z\in S^1, h(1,z)=1$,
\end{itemize}
realizes an homotopy between the identity of the circle and the constant function $z\mapsto 1\in S^1$.

More generally, the same kind of argument shows that every $n$-sphere is contractible. On the other hand, there exist a connected connectivity space $X$ such that no two distinct connectivity endomorphisms $X\rightarrow X$ are homotopic. For example, if $X=\mathcal{P}(\mathbf{R})$ is endowed with the integral connectivity structure for which non trivial connected subsets are subsets with a cardinal greater than the one of $\mathbf{R}$, then non-trivial connected subsets of $\mathbf{Cnct}(X,X)$ also have such a cardinal, and then every connectivity morphism from $\mathbf{I}$ to $\mathbf{Cnct}(X,X)$ is a constant function.

Those examples show that any theory of homotopy in the connectivity framework should be very different from the topological one. In particular, it could be interesting to use different kind of discrete times instead of $\mathbf{I}$.

\section{Pointed Connectivity Spaces}

\subsection{Pointed Sets} The category $\mathbf{pSet}$ of pointed sets and based maps is a concrete category on $\mathbf{Set}$. The forgetful functor $\mathbf{pSet}\to\mathbf{Set}$ will be denoted by $\vert-\vert$, and the base-point of a pointed set $P$ by $\beta(P)$, so $P=(\vert P\vert,\beta(P))$.

$\mathbf{pSet}$ has a zero object, $(\{*\},*)$, it is complete and cocomplete. In particular, the cartesian product of two pointed sets $P_1$ and $P_2$ is defined by $\vert P_1\times P_2\vert = \vert P_1\vert\times\vert P_2\vert$ and $\beta(P_1\times P_2)=(\beta(P_1),\beta(P_2))$. The class of coequalizers coincides with the class of all epimorphisms, \textit{i.e.} surjective based maps, and with the class of quotient morphisms (in $\mathbf{pSet}$ every morphism is final). If $\sim$ is an equivalence relation on $\vert P\vert$, the quotient pointed set $P/\sim$ is defined by $\vert P/\sim\vert=\vert P\vert/\sim$ and  $\beta(P/\sim)=\widetilde{\beta(P)}$. In particular, if $T$ is a subset of $\vert P\vert$, $P/T$ denotes the pointed set $P/\sim_T$. 
The coproduct of $P_1$ and $P_2$ is denoted by $P_1\vee P_2$. It can be defined either as the quotient of the set $\vert P_1\vert \amalg \vert P_2\vert$ by the equivalence relation which identifies $\beta(P_1)$ and $\beta(P_2)$ or alternatively by the formulas
\begin{equation} \label{smash}
\vert P_1\vee P_2\vert =  (\vert P_1\vert\times\{\beta(P_2)\})\cup (\{\beta(P_1)\}\times\vert P_2\vert)
\end{equation}
and
\[ 
\beta(P_1\vee P_2) =(\beta(P_1),\beta(P_2)).
 \]

The category $\mathbf{pSet}$ is not cartesian closed since, for example, if $P$ is a pointed set with two elements and $Q$ is the zero object, then $P\times (Q\vee Q)\simeq P$ whereas $(P\times Q)\vee(P\times Q)$ has three elements. Nevertheless, the set of based maps from a pointed set $P$ to a pointed set $Q$  has a ``natural''  special point, that is the constant map $x\mapsto \beta(Q)$, so there is a ``natural'' object in $\mathbf{pSet}$ representing $Hom(P,Q)$. Let \[\mathbf{pSet}(P,Q)=(Hom(P,Q),x\mapsto \beta(Q))\] denotes this object. For each pointed set $P$, we then have an endofunctor $\mathbf{pSet}(P,-)$ on $\mathbf{pSet}$, with $\mathbf{pSet}(P,f)=f\circ -$. One knows that this functor has a left adjoint $P\wedge -$, the so-called \textit{smash product}, defined on objects by 
\[P\wedge Q= (P\times Q)/\vert P\vee Q\vert,\] where the set $\vert P\vee Q\vert$ is defined by the formula (\ref{smash}), and on based maps $f:Q\to R$ by
\begin{equation} \label{smash function}
\forall (p,q)\in \vert P \vert \times \vert Q\vert, (P\wedge f)(\widetilde{(p,q)})= \widetilde{(p,f(q))}.
\end{equation}

Then, endowed with the smash product, $\mathbf{pSet}$ is a closed symmetric monoidal category. Note that there are no projections associated with the smash product, and that the two-elements pointed set is a unit for it.

\subsection{Pointed Integral Connectivity Spaces} 

\begin{df}
A \emph{pointed integral connectivity space} $X$ is a triple \linebreak $(S,\mathcal{K},b)$, where $(S,\mathcal{K})$ is an integral connectivity space and $b$ a point of $S$, called the base-point of $X$.
\end{df}

For every pointed connectivity space $X$, we will denote $\vert X\vert$ its underlying carrier set, $\kappa(X)$ its connectivity structure and $\beta(X)$ its base-point, so $X=(\vert X\vert,\kappa(X),\beta(X))$.

The category whose objects are the pointed integral connectivity spaces and whose morphisms are connectivity morphisms preserving base-points will be denoted by $\mathbf{pCnct}$. It can be viewed as a category with lattices of structures on the base category $\mathbf{pSet}$ of pointed sets. Indeed, the choice of a base-point does not have any effect on the lattice of (integral) connectivity structures on a given set, and connectivity morphisms between pointed spaces are just based maps between underlying pointed sets which preserve connected subsets, so  $\mathbf{pCnct}=\mathbf{pSet}_{pCnct}$ with $pCnct=Cnct\circ \vert - \vert :\mathbf{pSet}\to\mathbf{JCPos}$. Thus,

\begin{prop} $\mathbf{pCnct}$ is a topological category on $\mathbf{pSet}$. It is thus complete and cocomplete.
\end{prop}

The topological forgetful functor $\mathbf{pCnct}\to\mathbf{pSet}$ will be denoted $\vert - \vert_p$, so that $\vert X\vert_p=(\vert X\vert,\beta(X))$. The category $\mathbf{pCnct}$ can also be viewed as a concrete category on $\mathbf{Cnct}$, and we will denote $\vert-\vert_\kappa$ the corresponding forgetful functor, so that $\vert X \vert_\kappa=(\vert X\vert,\kappa(X))$.
Then, the product of two pointed integral connectivity spaces $X_1$ and $X_2$ is characterised by $\vert X_1\times X_2\vert_p = \vert X_1\vert_p\times \vert X_2\vert_p$ and $\vert X_1\times X_2\vert_\kappa = \vert X_1\vert_\kappa\times \vert X_2\vert_\kappa$. If $\sim$ is an equivalence relation on $\vert X\vert$, the quotient pointed space $X/_\sim$ is likewise characterised by $\vert X/_\sim\vert_p=\vert X\vert_p/_\sim$ and $\vert X/_\sim\vert_\kappa=\vert X\vert_\kappa/_\sim$. This gives in particular the definition of $X/T$ with $T\subseteq \vert X\vert$. The coproduct  satisfies $\vert X_1 \vee X_2\vert_p=\vert X_1\vert_p \vee \vert X_2\vert_p$, and its connectivity part $\vert X_1 \vee X_2\vert_\kappa$ can be defined either as the quotient of $\vert X_1\vert_\kappa\amalg\vert X_2\vert_\kappa$ by the relation $\beta(X_1)\sim \beta(X_2)$, or as induced by the space $\vert X_1\vert_\kappa \boxtimes \vert X_2\vert_\kappa$ on $\vert X_1 \vee X_2\vert$ seen as a subset of $\vert X_1\vert \times\vert X_2\vert$ according to the formula (\ref{smash}), the $X_i$ replacing there the $P_i$. In the sequel, the expression $\vert X_1\vee X_2\vert$ will keep this last meaning. Now, the same argument as for $\mathbf{pSet}$ shows that $\mathbf{pCnct}$ is not cartesian closed. 

\pagebreak[3]
\subsection{The Smash Product}
\begin{df} Let $X_1$ and $X_2$ be pointed integral connectivity spaces. Then,
\begin{itemize}
\item the \emph{tensor product}  $X_1\boxtimes X_2$ is defined by the relations 
\begin{enumerate}
\item $\vert X_1\boxtimes X_2\vert_p=\vert X_1\vert_p\times \vert X_2\vert_p$,
\item $\vert X_1\boxtimes X_2\vert_\kappa=\vert X_1\vert_\kappa\boxtimes \vert X_2\vert_\kappa$,
\end{enumerate}
\item the \emph{smash product} is defined by
$X_1\wedge X_2 = (X_1\boxtimes X_2) / \vert X_1 \vee X_2\vert$,
\item $\mathbf{pCnct}(X_1,X_2)$, the \emph{pointed connectivity space of connecting based maps} from $X_1$ to $X_2$, is defined by
\begin{enumerate}
\item $\vert\mathbf{pCnct}(X_1,X_2)\vert=\vert\mathbf{Cnct}(\vert X_1\vert_\kappa,\vert X_2\vert_\kappa)\vert\cap\vert\mathbf{pSet}(\vert X_1\vert_p,\vert X_2\vert_p)\vert$,
\item $\kappa(\mathbf{pCnct}(X_1,X_2))=i^*(\kappa(\mathbf{Cnct}(\vert X_1\vert_\kappa,\vert X_2\vert_\kappa)))$, where $i$ is the inclusion map $i:\vert\mathbf{pCnct}(X_1,X_2)\vert \hookrightarrow \vert\mathbf{Cnct}(\vert X_1\vert_\kappa,\vert X_2\vert_\kappa)\vert$,
\item $\beta(\mathbf{pCnct}(X_1,X_2))$ is the constant map $x\mapsto\beta(X_2)$.
\end{enumerate}
\end{itemize}
\end{df}

Now, with those objects we can define, for every pointed integral connectivity space $X$, the endofunctors $\mathbf{pCnct}(X,-)$ and $X\wedge -$ on the category $\mathbf{pCnct}$. In fact, for every morphism $f$, the morphisms $X\wedge f$ and $\mathbf{pCnct}(X,f)$ are given by the same formulas as for the corresponding endofunctors on $\mathbf{pSet}$. 

\begin {thm} For every pointed integral connectivity space $X$, the endofunctor $(X\wedge -)$ on $\mathbf{pCnct}$ is left adjoint to the endofunctor $\mathbf{pCnct}(X,-)$.
\end {thm}
\noindent \textit{Proof}. 
Let $X$, $Y$ and $Z$ be pointed integral connectivity spaces. For every based connecting map $\psi:X\wedge Y\to Z$, one has a based connecting map $\rho(\psi):Y\to \mathbf{pCnct}(X,Z)$ defined for all $y\in Y$ by 
\[ \rho(\psi)(y)=\psi(\widetilde{(-,y)}).  \]
Indeed, for every $y\in Y$, $\psi(\widetilde{(-,y)})\in\mathbf{pCnct}(X,Z)$ since
\begin{itemize}
\item $\psi$ is defined on classes $\widetilde{(x,y)}$, so $\psi(\widetilde{(-,y)})$ is a function from $\vert X\vert$ to $\vert Z\vert$,
\item  $\psi(\widetilde{(-,y)})(\beta(X))=\psi(\beta(X\wedge Y))=\beta(Z)$,
\item for every $K\in \kappa(X)$, $s(K\times\{y\})\in \kappa(X\wedge Y)$  so $\psi(\widetilde{(-,y)}(K)\in\kappa(Z)$,
\end{itemize}
where $s:X\boxtimes Y\twoheadrightarrow X\wedge Y$ denotes the canonical map.
And the function $y\mapsto\psi(\widetilde{(-,y)})$ is a based connecting map from $Y$ to $\mathbf{pCnct}(X,Z)$, since
\begin{itemize}
\item $\psi(\widetilde{(-,\beta(Y))})=(x\mapsto \beta(Z))=\beta(\mathbf{pCnct}(X,Z))$,
\item for every $L\in\kappa(Y)$, $\{\psi(\widetilde{(-,y)}),y\in L\}\in\kappa(\mathbf{pCnct}(X,Z))$, since for every $x\in\vert X\vert$ one has $<\{\psi(\widetilde{(-,y)}),y\in L\},x>=\psi(\widetilde{(x,L)})\in\kappa(Z)$.
\end{itemize}

\noindent Now, one verifies as well that the formula
\[ \theta(\varphi)(\widetilde{(x,y)})=\varphi(y)(x) \]
defines a map $\theta$ from $Hom(Y,\mathbf{pCnct}(X,Z))$ to $Hom(X\wedge Y,Z)$, and that $\theta$ and $\rho$ are inverses of each other.
Finally, $\rho$ is natural since for all pointed integral connectivity spaces $Y$, $Y'$, $Z$ and $Z'$ and for all based connecting maps $u:Y\to Y'$, $v:Z\to Z'$ and $\psi:X\wedge Y'\to Z$, one has $\rho(v\circ \psi\circ (X\wedge u))=(y\mapsto v\circ \psi(\widetilde{(-,u(y))}))=\mathbf{pCnct}(X,v)\circ \rho(\psi)\circ u$.
\begin{flushright}$\square$\end{flushright}

\section{Finite Integral Connectivity Spaces}

\subsection{Generic Graphs}

\begin{df} Let $X$ be a finite integral connectivity space. A \emph{generic point} of $X$ is a non-empty irreducible connected subset of $X$. The \emph{generic graph} $G_X$ of $X$ is the directed graph whose vertices are the generic points of $X$ and such that $g\rightarrow h$ is a directed edge of $G_X$ if and only if $g\supsetneqq h$ and there is no generic point $k$ such that  $g\supsetneqq k \supsetneqq h$.
\end{df}

Associated with a partial order, the directed graph $G_X$ is a so-called \emph{directed acyclic graph}, that is a directed graph with no \textit{directed} cycle; note that cycles are allowed in the undirected graph obtained by forgetting orientation of the edges. On the other hand, not every finite acyclic directed graph is a $G_X$ for some finite integral connectivity space $X$. For example, the directed acyclic graph $a\rightarrow b$ is not such a $G_X$.

\textit{Notation}. For the sake of simplicity, if $G$ is a directed graph,  $a\in G$ will express that $a$ is a vertex of $G$ and $(a\rightarrow b)\in G$ will express that $a\rightarrow b=(a,b)$ is a directed edge of this graph.

\begin{prop} A finite integral connectivity space $X$ is characterised, up to isomorphism, by its generic graph $G_X$ (defined up to isomorphism).
\end{prop}
\noindent \textit{Proof}. The space $X$ being integral, every singleton is an irreducible connected subset, and appears in $G_X$ as a sink, \textit{i.e.} a vertex with no outgoing edges. Thus, the carrier $\vert X\vert$ of the space is given, up to bijection, by the set of sinks of $G_X$. Now, the connectivity structure is given by $G_X$ as a consequence of the proposition \ref{characterization of finite spaces by irreducible subsets}. 
\begin{flushright}$\square$\end{flushright} 

\begin{prop} If $X$ is a non-empty finite integral connectivity space, then
\begin{enumerate}
\item $X$ is connected iff $G_X$ is connected,
\item there is a bijection between connected components of $X$ and those of $G_X$,
\item $X$ is irreducible iff $G_X$ has exactly one source, \textit{i.e.} a vertex with no incoming edges,
\item $X$ is distinguished iff there is no triple $(a,b,c)$ of distinct vertices in $G_X$ such that $(a\rightarrow b)$ and $(b\leftarrow c)$ are in $G_X$.
\item $X$ is connected and distinguished iff $G_X$ is a directed tree.
\end{enumerate}
\end{prop}
\noindent Proof. 
\begin{enumerate}
\item If there is an arrow $(a\rightarrow b)$ in $G_X$ then $a$ and $b$, as subsets of $\vert X\vert$, are containded in the same connected component of $X$; thus, if $G_X$ is connected then $X$ is also connected. 
On the other hand, let $(C_i)$ be the family of $G_X$ connected components and, for each $i$, let $\sigma(C_i)$ be the union of sinks belonging to $C_i$; then, every connected subset produced at any step of the process described in the theorem \ref{Generation of connectivity structures} stays in one of the $\sigma(C_i)$, otherwise there should be two irreducible connected subsets of $X$ contained respectively in two distinct $\sigma(C_i)$ and with a non-empty intersection, which is not possible. Thus, if $G_X$ is not connected, neither is $X$.

\item The generic graph $G_X$ of the disjoint union $X$ of any finite family of  finite spaces $X_i$ is clearly the disjoint union of the $G_{X_i}$, thus the connected components of any finite space $X$ are the $\sigma(C_i)$ associated with the connected components $C_i$ of $G_X$.

\item If $X$ is irreducible then $\vert X\vert$ is a generic point which contains all other generic points so it is the only source in $G_X$.

If $G_X$ has only one source, then each irreducible connected proper subset of $X$ is contained in a larger irreducible subset, so, $X$ being finite and the set of irreducible connected sets being nonempty, $\vert X\vert$ is itself an irreducible connected subset.

\item  If there is a triple $(a,b,c)$ with $a\neq c$ and $a\rightarrow b \leftarrow c$ in $G_X$, then $a \cup c$ is a reducible connected subset of $X$ which is thus not distinguished.

If two irreducible connected subsets of $X$ not included one in the other have a common point, then there must exist in $G_X$ a triple of distinct points $(a,b,c)$ with $a\rightarrow b \leftarrow c$ in $G_X$; thus, if $G_X$ does not admit such a triple, then the inductive generation of connected subsets from irreducible ones (theorem \ref{Generation of connectivity structures}) cannot produce any other connected set than the latters.
\item The last affirmation is a direct consequence of the others.
\end{enumerate}
\begin{flushright}$\square$\end{flushright}

\begin{df} Let $X$ be a non-empty finite  integral connectivity space. The \emph{index} of any irreducible subset of $X$ is its height as a vertex of the directed acyclic graph $G_X$ (i.e. the length of the longest path from that vertex to a sink of $G_X$). The index $\omega(X)$ of $X$ is the maximum of indexes of its irreducible connected subsets, that is the length of $G_X$.
\end{df}

\begin{example} A finite space of index $0$ is totally disconnected, \textit{i.e.} its structure is the discrete one.
\end{example}

\begin{example} One has $\omega(U_G(S))\leq 1$ for any finite simple undirected graph $S$.
\end{example}

The definition of the index of a finite integral connectivity space results in the definition of a new numerical invariant for links:

\begin{df} The \textit{connectivity index} of a tame link $L$ in $\mathbf{R}^3$ (or $\mathbf{S}^3$) is $\omega(L)=\omega(S_L)$.
\end{df}

\begin{example} The connectivity index of the Borromean link or, more generally, of any Brunnian link, is $\omega(\mathbf{B}_n)=1$.
\end{example}

\begin{rmq} The connectivity index is not a Vassiliev finite type invariant for links. For example, it is easy to check that the connectivity index of the \textit{singular} link with two components, a circle and another component crossing this circle at $2n$ double-points, is greater than $2^n$.
\end{rmq}

\begin{prop} 
One has $\omega(X)\leq \mathrm{card}(X)-1$ for every finite integral space $X$; and 
the integral connectivity space $\mathbf{V}_n$ defined by $\vert \mathbf{V}_n\vert=n$ and $\kappa(\mathbf{V}_n)^\bullet=\{2,3,\cdots,n\}$ is, 
up to isomorphism, the only  integral connectivity space such that $\mathrm{card}(\mathbf{V}_n)=n$  and $\omega(\mathbf{V}_n)=n-1$.
\end{prop}
\noindent \textit{Proof}. 
A trivial induction results in the first claim.
The second one is obvious if $n=1$. Suppose that it is true for an integer $n$, and let $X$ be an integral connectivity space with $n+1$ points and with index $n$. Then there must exist an irreducible connected subset $K$ of $X$ with index $n-1$, and one has necessarily $\mathrm{card}(K)\geq n$, so $\mathrm{card}(K)=n$. By induction, $K\backsimeq \mathbf{V}_n$. Let $x$ be the unique element of $X\setminus K$. $\vert X\vert$ is necessarily the only non-trivial connected subset which contains $x$, otherwise $X$ would be of index smaller than $n$, then  $\kappa(X)=\{\{x\}\}\cup\kappa(K)\cup \{\vert X\vert\}$, and thus $X\simeq \mathbf{V}_{n+1}$.
\begin{flushright}$\square$\end{flushright} 

Let us now describe two ways to product new finite spaces from two given non-empty finite integral connectivity spaces $X$ and $Y$, $Y$ being supposed irreducible.
\begin{enumerate}
\item Let $x$ be a point of $\vert X\vert$. We denote by $X\rhd_x Y$ the connectivity space whose generic graph is obtained by replacing in $G_X$ the sink $\{x\}$ by (a copy of) $G_Y$, arrows to $x$ in $G_X$ being replaced by arrows to the unique source of (the copy of) $G_Y$. In other words, $X\rhd_x Y$ is the integral space such that $\vert X\rhd_x Y\vert=\vert X\vert\smallsetminus\{x\}\cup \vert Y'\vert$ and the set $\kappa_0(X\rhd_x Y)$ of irreducible connected sets is given by \[\{K\in\kappa_0(X), x\notin K\}\cup\kappa_0(Y')\cup \{K\cup \vert Y'\vert, x\in K\in\kappa_0(X) \},\] where  $Y'$ is a copy of $Y$ such that $\vert X\vert\cap \vert Y'\vert=\emptyset$.
\item We can replace simultaneously every sink of $G_X$ by (a copy of) $G_Y$ to produce a space denoted by $X\rhd Y$. That is, $X\rhd Y$ is the connectivity space such that $\vert X\rhd Y\vert = \vert X\vert \times \vert Y\vert$ and the set $\kappa_0(X\rhd Y)$ of irreducible connected sets is given by
\[\kappa_0(X\rhd Y)=\{\{x\}\times L, x\in \vert X\vert, L\in\kappa_0(Y)\}\cup\{K\times \vert Y\vert, K\in\kappa_0(X)\}.\]

\end{enumerate}

\begin{example} $\mathbf{B}_2\rhd_x \mathbf{V}_n\simeq \mathbf{V}_{n+1}$, where $x$ is any of the two points of $\mathbf{B}_2$.
\end{example}

\begin{prop} For any non-empty finite integral connectivity space $X$ and any non-empty irreducible finite integral connectivity space $Y$, one has $\omega(X\rhd Y)=\omega(X)+\omega(Y)$.
\end{prop}
\noindent \textit{Proof}. By construction, $G_{X\rhd Y}$ is obtained by replacing each sink of $G_X$ by a copy of $G_Y$, so its length is  $\omega(X)+\omega(Y)$.
\begin{flushright}$\square$\end{flushright}

\begin{example}
The link depicted on figure \ref{borroborro} is a Borromean assembly of three Borromean links. Its generic graph is (isomorphic to) $\mathbf{B}_3\vartriangleright \mathbf{B}_3$, and its connectivity index is $2$.
\begin{figure} 
\begin{center}
\includegraphics [scale=0.2]{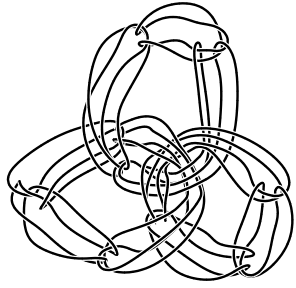}
\caption{A Borromean ring of borromean rings.}
\label{borroborro}
\end{center}
\end{figure}
\end{example}

\subsection{Representation by Links}

In \cite{Dugowson:2007c, Dugowson:20070717}, I asked whether every finite connectivity space can be represented by a link, \textit{i.e.} whether there exists a link whose connectivity structure is (isomorphic to) the one given. It turns out that in 1892, Brunn \cite{Brunn:1892a} first asked this question, without clearly bringing out the notion of a connectivity space. His answer was positive, and he gave the idea of a proof based on a construction using some of the links now called ``Brunnian''. In 1964, Debrunner \cite{Debrunner:1964}, rejecting the Brunn's ``proof", gave another construction, proving it but only for $n$-dimensional links with $n\geq 2$. In 1985, Kanenobu \cite{Kanenobu:198504,Kanenobu:1986} seems to be the first to give a proof of the possibility of representing every finite connectivity structure by a classical link, a result which is still little known at this date. The key idea of those different constructions is already in the Brunn's original article; it consists in using some Brunnian structures to successively link the sets of components which are desired to become unsplittable. Thus, in Brunn's point of view, the links called today ``Brunnian links" are not so interesting in themselves, but more for the constructions they allow to make, that is the representation of \textit{all} finite connectivity strutures by links. 

\begin{thm} [Brunn-Debrunner-Kanenobu] Every finite connectivity structure is the splittability structure of at least one link in $\mathbf{R}^3$.
\end{thm}

\begin{rmq} Note that the structure of the links used by Brunn is well described by the so-called \textit{Brunnian groups} constituted by the \textit{Brunnian braids} introduced as \textit{decomposable braids} by Levinson \cite{Levinson:1973,Levinson:1975} 
(see also \cite{Stanford:1999} and \cite{Wu:2007})
and by the \textit{Brunnian words} studied by Gartside and Greenwood \cite{Gartside.Greenwood:20031205,Gartside.Greenwood:20060606}.  \end{rmq}

\begin{figure} 
\begin{center}
\includegraphics [scale=0.3]{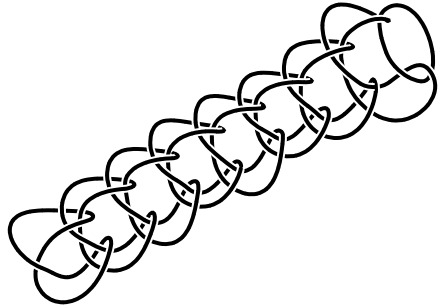}
\caption{A link with a connectivity index 8.}
\label{tendu}
\end{center}
\end{figure}

\begin{example} The structure of the connectivity space $\mathbf{V}_9$ with $9$ points and maximal connectivity index $8$ is the splittability structure of the link depicted on figure \ref{tendu}.
\end{example}

\noindent \textbf{Acknowledgments}. Thanks to Jean Bénabou, who introduced me to categories with lattices of structures. To Sergei Soloviev, who asked me whether smash products were possible for connectivity spaces. To David C. Ullrich who, on the forum \texttt{sci.math}, gave the upper bound $\omega_0$ for the construction I give in the theorem \ref{Generation of connectivity structures}. To René Guitart, Mark Weber, Albert Burroni and Quentin Donner for various talks. To Christopher-David Booth, Behrouz Roumizadeh and Anne Richards who helped me to correct my English.


\tableofcontents

\end{document}